\documentclass[11pt]{article}
\usepackage[T1]{fontenc}
\usepackage{graphicx}%
\usepackage{amsmath,amssymb,amsfonts}%
\usepackage{mathrsfs}%
\usepackage[title]{appendix}%
\usepackage{textcomp}%
\usepackage{manyfoot}%
\usepackage{booktabs}%
\usepackage{algorithm}%
\usepackage{algorithmicx}%
\usepackage{algpseudocode}%
\usepackage{listings}%
\usepackage[table,xcdraw]{xcolor}
\usepackage{caption}
\usepackage{subcaption}
\usepackage{float}
\usepackage{mathptmx}
\usepackage{lmodern}
\usepackage{wrapfig}

\date{}
\begin{document}
	\title{Reconstruction of Boundary Data in the Helmholtz Equation Using Particle Swarm Optimization}
	
	\author{Jamal Daoudi$^1$ and Chakir Tajani$^{1,*}$\\
		$^1$Department of Mathematics, Polydisciplinary faculty of Larache,\\
		 Abdelmalek Essadi University, Morocco\\
         jamaldaouid5@gmail.com \\
		ch.tajani@uae.ac.ma}
	
	\maketitle

\begin{abstract}
This paper tackles the data completion problem related to the Helmholtz equation. The goal is to identify unknown boundary conditions on parts of the boundary that cannot be accessed directly, by making use of measurements collected from accessible regions. Such inverse problems are known to be ill-posed in the Hadamard sense, which makes finding stable and dependable solutions particularly difficult. To address these challenges, we propose a bio-inspired method that combines Particle Swarm Optimization with Tikhonov regularization. The results of our numerical experiments suggest that this approach can yield solutions that are both accurate and stable, converging reliably. Overall, this method provides a promising way to handle the inherent instability and sensitivity of these types of inverse problems.



\end{abstract}

\maketitle

\section{Introduction}
\label{intro}
The Helmholtz equation is fundamental in describing wave propagation and vibration phenomena, with applications across a wide range of scientific and engineering disciplines. It commonly appears in models of wave scattering, electromagnetic scattering, and structural vibrations, including problems involving acoustic cavities and related scattering scenarios. Because of its broad importance, the Helmholtz equation has attracted considerable attention in research, leading to the development of numerous solution methods. In particular, the Cauchy problem for the Helmholtz equation plays a key role in advancing our understanding of complex wave behaviors and continues to be actively investigated \cite{beskos1997boundary,chen1998dual,hall1995boundary,hong2018three,marin2003conjugate,wood1995steady,harari1998boundary}.

Despite more than a century of extensive research on the Helmholtz equation, challenges remain when boundary data are incomplete or inaccessible. In such cases, inverse problems arise, often complicated by noise and partial data, which make accurate determination of boundary conditions difficult. The problem under consideration exemplifies an ill-posed inverse problem: small perturbations in the input data can cause large deviations in the solution \cite{hadamard1923lectures}. Due to the inherent ill-conditioning of the Cauchy problem for the Helmholtz equation, conventional numerical methods often fail, requiring the use of specialized regularization techniques. Over the years, various regularization methods have been proposed to address this challenge, including the Landweber method \cite{marin2004bem}, conjugate gradient methods \cite{marin2003conjugate,marin2004comparison}, Fourier regularization \cite{elden2000wavelet,fu2009fourier}, boundary knot methods \cite{jin2005boundary}, plane wave methods \cite{jin2008plane}, spherical wave expansions \cite{yu2008acoustic}, boundary element-minimal error methods \cite{marin2009boundary}, and fundamental solution methods \cite{marin2005method}. These approaches typically focus on finding deterministic solutions that best fit the given data within an optimization framework.\\

In contrast, metaheuristic methods adopt probabilistic search strategies to explore a wider solution space, identifying optimal candidates from many possibilities. Such techniques are especially valuable for ill-posed inverse problems because they can accommodate multiple potential solutions consistent with the available data. Examples include artificial bee colony optimization \cite{karaboga2014comprehensive}, ant colony optimization \cite{de1988learning}, particle swarm optimization \cite{kennedy1995particle}, the Bat algorithm \cite{socha2008ant}, and genetic algorithms \cite{yang2012bat}. These algorithms are notable for their flexibility in handling diverse inverse problems, including those with non-smooth or non-convex objective functions, as well as their robustness in the presence of noisy or incomplete data. However, achieving optimal performance is highly dependent on careful tuning of the algorithm parameters and the appropriate selection of methods.\\

This study develops a computational algorithm for solving the Helmholtz equation by integrating particle swarm optimization (PSO) with the least squares method, while employing the finite element method (FEM) to solve the direct problem. The proposed approach uses the global search capability of PSO along with the precision of FEM, using Cauchy data from the accessible boundary to improve the precision and efficiency of the inverse solution.\\

The remainder of the paper is organized as follows. Section \ref{sec:2} details the mathematical formulation of the inverse problem and its reformulation as an optimization problem. Section \ref{sec:3} explains the application of PSO to solve this problem. Section \ref{sec:4} evaluates the stability, accuracy, and efficiency of the proposed method through numerical experiments in both regular and irregular domains. Finally, Section \ref{sec:5} presents the main conclusions and discusses future directions.

\section{Mathematical Formulation}
\label{sec:2}

Let $\Omega$ denote an open and bounded domain in $\mathbb{R}^d$ $(d=2,3)$ with a smooth boundary $\Gamma$. The boundary is divided into two disjoint parts, $\Gamma = \Gamma_i \cup \Gamma_c$, where $\Gamma_i \cap \Gamma_c = \emptyset$ and $\operatorname{mes}(\Gamma_c) \neq 0$.

The data completion problem for the Helmholtz equation is mathematically formulated as:
\begin{equation}
(P): \left\{
\begin{aligned}
-\Delta u - \kappa^2 u = 0 & & & \text{ in } \Omega, \\ 
\partial_n u = g & & & \text{ on } \Gamma_c, \\
u = f & & & \text{ on } \Gamma_c,
\end{aligned}
\right.
\label{Problem}
\end{equation}
where $\Delta$ is the Laplacian operator, $n$ denotes the outward unit normal to $\Gamma$, $\kappa$ is the wave number (a complex value), and $\kappa = \mathrm{i}\lambda$ (with $\mathrm{i} = \sqrt{-1}$) in the case of the modified Helmholtz equation. The terms $\partial_n$ represent the outward normal derivative and $f$ and $g$ denote the Cauchy data available on the accessible boundary $\Gamma_c$.\\

\textbf{Direct Problem:} \\
The direct problem aims to find the solution \( u \) by solving the Cauchy problem given in Eq.~\eqref{Problem}, utilizing the prescribed Cauchy data—either \( u \) and its normal derivative \( \partial_n u \) on the accessible boundary segment \( \Gamma_c \) and inaccessible boundary segment \( \Gamma_i \), respectively, or vice versa.

\vspace{0.3cm}

\textbf{Inverse Problem:} \\
The inverse problem involves recovering the unknown values of \( u \) and its normal derivative on the inaccessible boundary \( \Gamma_i \) from the known data \( f \) and \( g \) specified on the accessible boundary \( \Gamma_c \).

\section{Optimization Problem}
\label{subsec:3}

Since the boundary value \(\phi_D = u|_{\Gamma_i}\) on the inaccessible boundary \(\Gamma_i\) is unknown, we must first solve the following direct problem:

\begin{equation}
(P_D): 
\begin{cases}
-\Delta u - \kappa^2 u = 0 & \text{in } \Omega, \\ 
u = \phi_D & \text{on } \Gamma_i, \\
\partial_n u = g & \text{on } \Gamma_c,
\end{cases}
\label{PD}
\end{equation}

where, for given data \(\phi_D \in H^{1/2}(\Gamma_i)\) and \(g \in H^{-1/2}(\Gamma_c)\), there exists a unique solution \(u(\phi_D, g)\) to the problem \eqref{PD} (see \cite{evans2010vol}).

The main goal is to determine the unknown boundary value \(\phi_D\) such that the solution satisfies the following conditions on the accessible boundary \(\Gamma_c\):

\begin{equation}
\begin{cases}
u(\phi_D, g) = f & \text{on } \Gamma_c, \\
\partial_n u(\phi_D, g) = g & \text{on } \Gamma_c.
\end{cases}
\label{J}
\end{equation}

This problem can be reformulated as the minimization of the regularized least-squares functional \(\mathcal{J}_{DR}\), defined by

\begin{equation}
\mathcal{J}_{DR}(\phi_D) = \frac{1}{2} \|u(\phi_D, g) - f\|_{L^2(\Gamma_c)}^2 + \frac{\alpha}{2} \|\phi_D\|_{H^{1/2}(\Gamma_i)}^2,
\label{JDR}
\end{equation}

where the parameter \(\alpha > 0\) plays a critical role in balancing the fidelity to the measured data and the regularity of the solution. Specifically, the second term acts as a Tikhonov regularization that penalizes irregular or unstable solutions by imposing smoothness through the Sobolev norm \(\|\cdot\|_{H^{1/2}(\Gamma_i)}\). This regularization is essential to enhance the stability of the minimization process and to ensure well-posedness of the inverse problem.

\section{Solving the Cauchy Problem for the Helmholtz Equation Using a Stochastic Approach}
\label{sec:3}

\subsection{Particle Swarm Optimization}

Particle Swarm Optimization (PSO) is a heuristic algorithm inspired by the social behavior observed in natural swarms, such as bird flocks or fish schools. Introduced by Eberhart and Kennedy in 1995 \cite{kennedy1995particle}, the method relies on a population of candidate solutions—called particles—that explore the search space collectively. Each particle updates its position by considering both its own experience and the best solution found by the entire swarm, which guides the search towards optimal solutions.

Formally, for each particle \(i\) in an \(n\)-dimensional space, its velocity and position at iteration \(t\) are updated according to the equations:

\begin{equation}
\boldsymbol{V}_i(t+1) = \omega \boldsymbol{V}_i(t) + c_1 r_1 \left(\boldsymbol{p}_i - \boldsymbol{X}_i(t)\right) + c_2 r_2 \left(\boldsymbol{p}_g - \boldsymbol{X}_i(t)\right),
\label{EqP1}
\end{equation}

\begin{equation}
\boldsymbol{X}_i(t+1) = \boldsymbol{X}_i(t) + \boldsymbol{V}_i(t+1),
\label{EqP2}
\end{equation}

where \(\boldsymbol{X}_i(t)\) and \(\boldsymbol{V}_i(t)\) represent the position and velocity vectors of the particle at iteration \(t\). The terms \(\boldsymbol{p}_i\) and \(\boldsymbol{p}_g\) denote the best position found so far by particle \(i\) and by the swarm, respectively. The coefficients \(c_1\) and \(c_2\) control the influence of these best positions on the particle's movement, while the inertia weight \(\omega\) helps balance between exploring new regions and exploiting known good areas of the search space. Random variables \(r_1\) and \(r_2\), drawn uniformly from \([0,1]\), introduce variability to avoid premature convergence and enhance search diversity.

The PSO algorithm begins with a population of particles randomly distributed in the search space. At each iteration, the fitness of each particle is evaluated according to the objective function to be optimized. Based on these fitness values, each particle's personal best position and the global best position among the swarm are updated. Particle velocities and positions are adjusted according to the update equations. The process continues until a stopping criterion, such as a maximum number of iterations or minimum change in fitness, is met (as seen in Algorithm \ref{alg:pso}). PSO has been successfully applied to a variety of inverse problems, such as image reconstruction \cite{dash2012particle}, electrical impedance tomography \cite{wang2020image}, inverse scattering problems \cite{yang2021hybrid}, inverse transient radiation problems \cite{qi2011inverse}, and seismic inversion \cite{zhang2018ava, wu2019improved}.

\begin{algorithm}[H]
\caption{Particle Swarm Optimization (PSO)}
\label{alg:pso}
\begin{algorithmic}[1]
\Require $J(x)$ (objective function), $N$ (number of particles), $c_1, c_2$ (acceleration coefficients), $\omega$ (inertia coefficient), $maxiter$ (maximum iterations), $tolerance$ (convergence tolerance)
\Ensure $x^*$ (optimal solution), $\min J(x^*)$ (minimum value of objective function)
\State Initialize particle velocities $v_i$ and positions $x_i$ for $i = 1, \ldots, N$
\State Initialize the best positions $p_i = x_i$ and global best position of the swarm $p_g$
\State Evaluate $J(x_i)$ for each particle
\State Update $p_i$ and $p_g$
\For{$k = 1$ to $maxiter$}
    \For{$i = 1$ to $N$}
        \State Update velocity: $v_i = \omega v_i + c_1 r_1 (p_i - x_i) + c_2 r_2 (p_g - x_i)$
        \State Update position: $x_i = x_i + v_i$
        \State Evaluate $J(x_i)$
        \If{$J(x_i) < J(p_i)$}
            \State Update $p_i$
            \If{$J(p_i) < J(p_g)$}
                \State Update $p_g$
            \EndIf
        \EndIf
    \EndFor
    \If{$|J(p_g) - J(p_{{\text{old}}_g})| < tolerance$}
        \State \textbf{break}
    \EndIf
\EndFor
\State \Return $p_g$, $J(p_g)$
\end{algorithmic}
\end{algorithm}

\subsection{Algorithm for Solving the Cauchy Problem for the Helmholtz Equation}

In this section, we detail the steps of the PSO algorithm applied to the considered inverse problem. Specifically, given the Cauchy data \(g\) and \(f\) on the accessible boundary \(\Gamma_c\), the aim is to reconstruct the trace \(\phi_D = u_{|\Gamma_i}\) and the normal derivative \(\phi_N = \partial_n u_{|\Gamma_i}\) using the PSO algorithm in conjunction with the Finite Element Method (FEM). The PSO algorithm iteratively updates particle positions and velocities in the search space, evaluating the updated positions with FEM to obtain corresponding fitness values.

\begin{algorithm}[H]
\caption{Solving the Cauchy Problem for the Helmholtz Equation Using PSO}
\begin{algorithmic}[1]
\State \textbf{Step 1:} Parameter setting:
\State \hspace{0.5cm} \(N\): Swarm size
\State \hspace{0.5cm} \(c_1\) and \(c_2\): Acceleration coefficients
\State \hspace{0.5cm} \(\omega\): Inertia coefficient
\State \hspace{0.5cm} \(MaxIt\): Maximum number of iterations
\State \hspace{0.5cm} \(Lb\): Lower bound of the search space
\State \hspace{0.5cm} \(Ub\): Upper bound of the search space

\State \textbf{Step 2:} Initialization: 
\State \hspace{0.5cm} Randomly initialize the starting positions \(\phi_{D_k^{(0)}}\) for \(k = 1, \dots, N\) for each particle.
\State \hspace{0.5cm} Initialize each particle's velocity as \(v_k(0) = 0\).
\State \hspace{0.5cm} Set the best position of each particle \(p_k(0) = \phi_{D_k^{(0)}}\).

\State \textbf{Step 3:} Solve the direct problem \((P_D)\) for each given \(\phi_{D_k^{(0)}}\) using FEM:
\State \hspace{0.5cm} \begin{equation*}
(P_D): 
\begin{cases}
-\Delta u - \kappa^2 u = 0 & \text{in} \quad \Omega, \\ 
u = \phi_{D_k^{(0)}} & \text{on} \ \Gamma_i, \\
\partial_n u = g & \text{on} \ \Gamma_c.
\end{cases}
\end{equation*}
\State \textbf{Step 4:} Evaluation: 
\State \hspace{0.5cm} Compute the value of the objective function using Equation \eqref{JDR} for each particle \(\phi_{D_k^{(0)}}\).
\State \hspace{0.5cm} Set \(p_g^{(0)}\) to the particle with the lowest objective value in the swarm: 
\State \hspace{1cm} \(p_g(0) = \operatorname{argmin} \mathcal{J_{DR}}(\phi_{D_k^{(0)}})\).

\State \textbf{Step 5:} Update the positions and velocities of the particles:
\For{\(t = 1\) to \(MaxIt\)}
    \For{each particle \(k = 1\) to \(N\)}
        \State Generate two random numbers \(r_1, r_2\) between 0 and 1.
        \State Update the velocity using Eq. \eqref{EqP1}.
        \State Update the position using Eq. \eqref{EqP2}.
        \If{\(\mathcal{J_{DR}}(\phi_{D_{k^{(t)}}}) < \mathcal{J_{DR}}(p_k^{(t)})\)}
            \State Update \(p_{k^{(t)}} = \phi_{D_k^{(t)}}\).
        \EndIf
        \If{\(\mathcal{J_{DR}}(\phi_{D_{k^{(t)}}}) < \mathcal{J_{DR}}(p_{g^{(t)}})\)}
            \State Update \(p_{g^{(t)}} = p_k^{(t)}\).
        \EndIf
    \EndFor
\EndFor
\State \textbf{Step 6:} Output \(p_g\), the best solution found.
\end{algorithmic}
\label{algo2}
\end{algorithm}

\vspace{0.5em}
\noindent\textbf{Extension to the Neumann Condition.}  
This procedure focuses on reconstructing the Dirichlet boundary condition \(\phi_D\) on \(\Gamma_i\). To instead reconstruct the Neumann condition \(\phi_N = \partial_n u|_{\Gamma_i}\), only minor modifications are required:

\begin{itemize}
    \item In Steps 1 and 2, replace \(\phi_D\) with \(\phi_N\)
    \item In Step 3, replace \((P_D)\) with the Neumann problem \((P_N)\):
\[
(P_N): 
\begin{cases}
-\Delta u - \kappa^2 u = 0 & \text{in } \Omega, \\ 
\partial_n u = \phi_{N_k}^{(0)} & \text{on } \Gamma_i, \\
u = f & \text{on } \Gamma_c.
\end{cases}
\]
    \item In Step 4, use the following cost functional:
\begin{equation}
\mathcal{J}_{NR}(\phi_N) = \frac{1}{2} \|u(\phi_N, f) - g\|^2_{L^2(\Gamma_c)} + \frac{\alpha}{2} \|\phi_N\|^2_{H^{1/2}(\Gamma_i)}
\label{JNR}
\end{equation}
\end{itemize}

This adjusted formulation allows the PSO-FEM framework to handle both Dirichlet and Neumann data recovery for the Cauchy problem in the Helmholtz equation.

\newpage

\section{Numerical Results and Discussion}
\label{sec:4}

The goal of this work is to estimate Neumann and Dirichlet boundary conditions in situations where the analytical solution is not explicitly known. To this end, we employ a polynomial-based approximation method. The robustness and convergence of the proposed approach are validated through two two-dimensional test scenarios.

\subsection*{Case 1: Square Domain}

\begin{figure}[H]%
    \centering
    \subfloat[\centering Unstructured mesh over unit square.]{{\includegraphics[width=6cm]{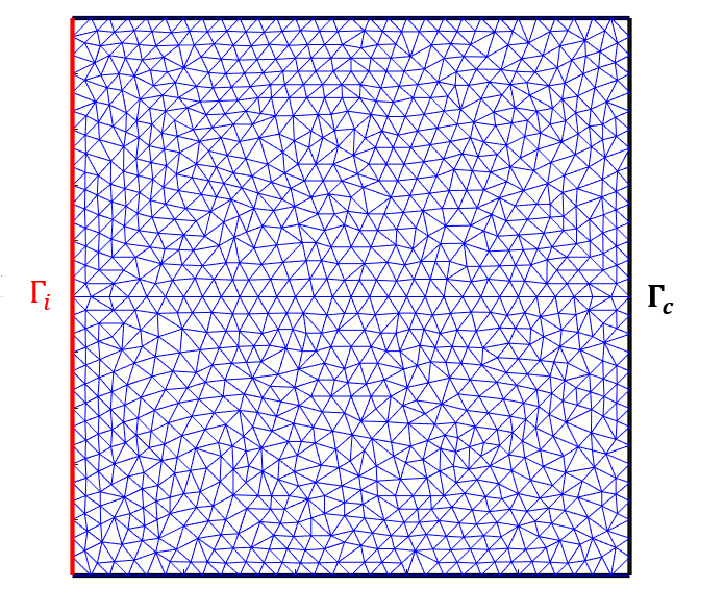} }\label{Square}}%
    \quad
    \subfloat[\centering Exact solution representation.]{{\includegraphics[width=6cm]{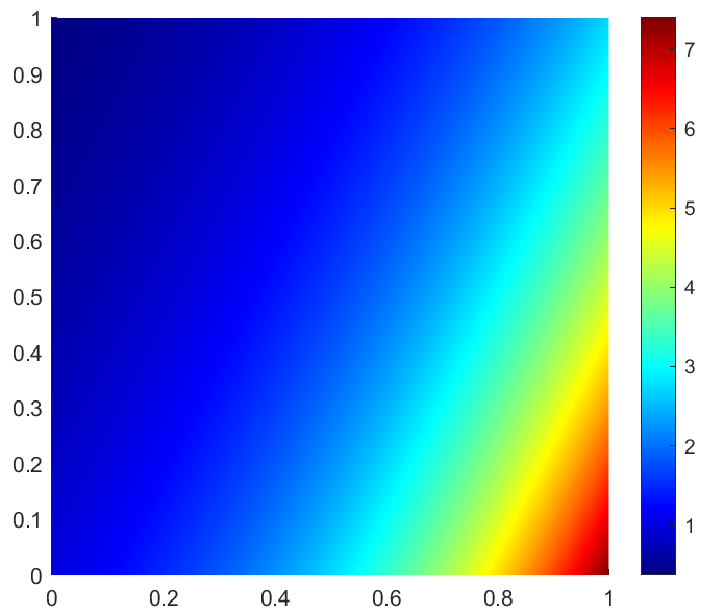} }\label{figquare2}}%
    \caption{Configuration of the Unit Square Domain}%
    \label{Solsquare}
\end{figure}

The first test is performed on a unit square domain, denoted by $\Omega = (0,1)^2$ (see Fig.~\ref{Square}). The boundary $\Gamma$ is split into two subsets: the inner boundary $\Gamma_i = \{(0,y): 0 < y < 1\}$ and the complementary part $\Gamma_c = \Gamma \setminus \Gamma_i$.

\subsection*{Case 2: Circular Domain}

\begin{figure}[H]%
    \centering
    \subfloat[\centering Unstructured mesh over unit disc.]{{\includegraphics[width=6cm]{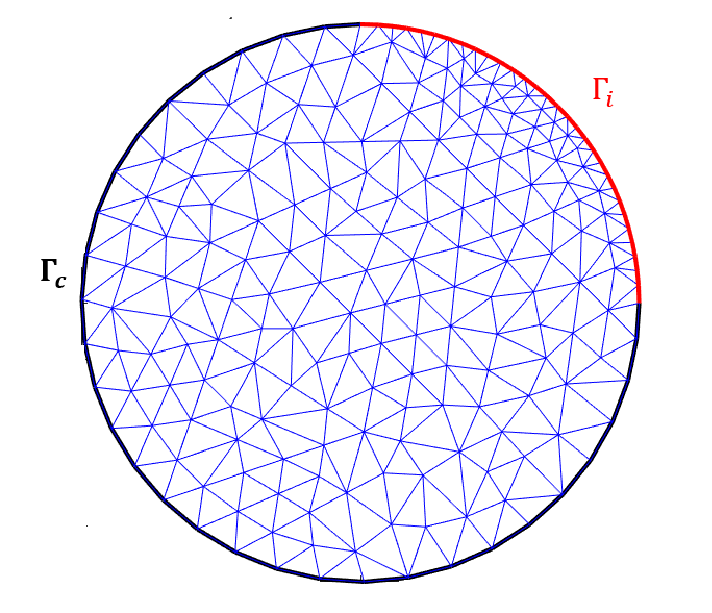} }\label{Cercle}}%
    \quad
    \subfloat[\centering Exact solution representation.]{{\includegraphics[width=6cm]{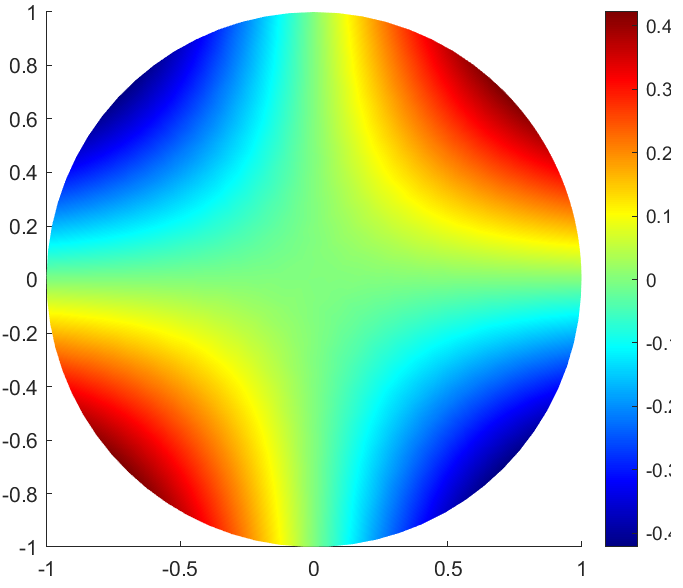} }\label{figCercle2}}%
    \caption{Configuration of the Unit Disc Domain}%
    \label{SolCercle}
\end{figure}

The second example involves a unit disc (Fig.~\ref{Cercle}). Here, the boundary is divided as follows:  
$\Gamma_i = \left\{(x(r,\theta), y(r,\theta)): x^2 + y^2 = r^2, \, \theta \in [0, \pi/2]\right\}$, and the remaining portion is assigned to $\Gamma_c = \Gamma \setminus \Gamma_i$.

In both configurations, we assume that the length (or measure) of $\Gamma_c$ is not less than that of $\Gamma_i$. The Dirichlet and Neumann boundary conditions on $\Gamma_i$ are derived analytically based on the exact solution used in each case:

\begin{equation}
\textbf{Case 1: } \left\{
\begin{array}{lll}
    u_{ex}(x,y) &= \exp(2x - y), & \text{with } k^2 = 5 \quad \text{(Modified Helmholtz)} \\
    \phi_D &= u_{ex}(x,y) \big|_{\Gamma_i}, & \\
    \phi_N &= \partial_n u_{ex}(x,y) \big|_{\Gamma_i}, &
\end{array}
\right.
\label{Ben1}
\end{equation}

\begin{equation}
\textbf{Case 2: } \left\{
\begin{array}{lll}
    u_{ex}(x,y) &= \sin(x)\sin(y), & \text{with } k^2 = -2 \quad \text{(Helmholtz)} \\
    \phi_D &= u_{ex}(x,y) \big|_{\Gamma_i}, & \\
    \phi_N &= \partial_n u_{ex}(x,y) \big|_{\Gamma_i}, &
\end{array}
\right.
\label{Ben2}
\end{equation}

Figures~\ref{figquare2} and~\ref{figCercle2} visualize the exact solutions within the entire computational domains.

\subsection*{Computational Setup}

All simulations were executed on a system equipped with an Intel(R) Core(TM) i7-8565U processor (1.80–1.99 GHz) and 16 GB of RAM. The implementation was carried out using the FreeFem++ platform~\cite{hecht2012new}.

\subsection*{PSO Configuration}

The parameters used for the Particle Swarm Optimization (PSO) algorithm are detailed below:

\begin{itemize}
    \item Swarm size: $N = 60$
    \item Acceleration constants: $c_1 = 1.5$, $c_2 = 1.5$
    \item Inertia weight: $\omega = 0.5$
    \item Maximum iterations: $MaxIt = 200$
    \item Search bounds: $Lb = -7$, $Ub = 7$
\end{itemize}

\subsection*{Synthetic Noise in Boundary Data}

In practical inverse problems, boundary measurements often contain errors due to physical or instrumental limitations. To simulate such uncertainty, we perturb the boundary data on $\Gamma_c$ using a noise model:

\begin{equation}
    T_{\text{noise}}(y) = T(y) \times \left(1 + \theta \cdot \nu\right), \quad \text{for} \quad y \in \Gamma_c
\end{equation}

Here, $\theta$ is a uniformly distributed random variable in the interval $[-1,1]$, and $\nu$ represents the noise intensity. In our implementation, random noise is generated using the FreeFem++ built-in function \texttt{randreal1()}.

\subsection{Unit Square}
\subsubsection{Choice of Regularization Parameter}

\begin{table}[h!] 
\centering
\begin{tabular}{|c|c|c|c|c|c|c|}
\hline
$\eta_{\alpha,\beta}$ & 1e-03       & 1e-04       & 1e-05       & 1e-06       & 1e-07       & {\color{red} 1e-08}       \\ \hline
$\mathcal{J_{DR}}(\phi_D)$ & 3.8687e-05 & 1.1567e-05 & 5.5405e-05 & 3.5806e-05 & 1.5557e-05 & {\color{red} 1.291e-05} \\ \hline
$\mathcal{J_{NR}}(\phi_N)$ & 6.6045e-04  & 5.1047e-04  & 5.1966e-04  & 5.2013e-04  & 5.4883e-04  & {\color{red} 4.9387e-04} \\ \hline
\end{tabular}
\caption{Values of $\mathcal{J_{DR}}(\phi_D)$ and $\mathcal{J_{NR}}(\phi_N)$ for different regularization parameters $\eta$.}
\label{RegSquare}
\end{table}

Table \ref{RegSquare} displays the values of the cost functions $\mathcal{J}_{DR}(\phi_D)$ and $\mathcal{J}_{NR}(\phi_N)$ for various regularization parameters $\eta_{\alpha,\beta}$, showing how sensitive the reconstruction is to this parameter. As $\eta$ decreases from $10^{-3}$ to $10^{-8}$, both cost functions generally decline, indicating improved alignment with the observed data. However, the rate of improvement slows for smaller values of $\eta$, suggesting a trade-off between data fidelity and solution smoothness. The optimal results are obtained at $\eta_{\alpha,\beta} = 10^{-8}$, where $\mathcal{J}_{DR}(\phi_D)$ reaches $1.291 \times 10^{-5}$ and $\mathcal{J}_{NR}(\phi_N)$ achieves $4.9387 \times 10^{-4}$, marking the most accurate and stable reconstructions across the tested range. These findings confirm that only a minimal amount of regularization is needed to achieve reliable results while avoiding overfitting. Therefore, the values $\eta_{\alpha,\beta} = 10^{-8}$ are adopted in the following numerical experiments for both Dirichlet and Neumann formulations to ensure consistency and optimal performance.

\subsubsection{Convergence of the Algorithm}
\begin{figure}[H]
    \centering
    \begin{subfigure}[b]{0.45\textwidth}
        \includegraphics[width=\textwidth]{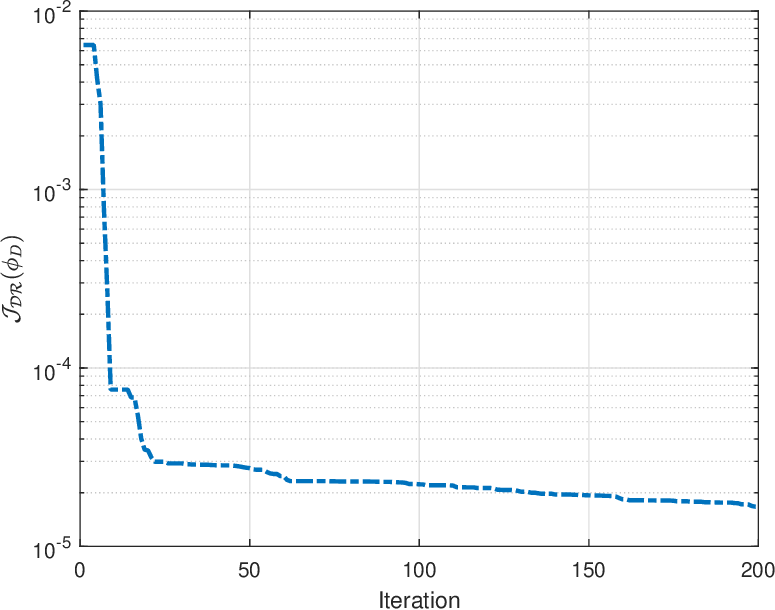}
        \caption{$\mathcal{J_{DR}}(\phi_D)$.}
        \label{fig1}
    \end{subfigure}
    \hfill
    \begin{subfigure}[b]{0.45\textwidth}
        \includegraphics[width=\textwidth]{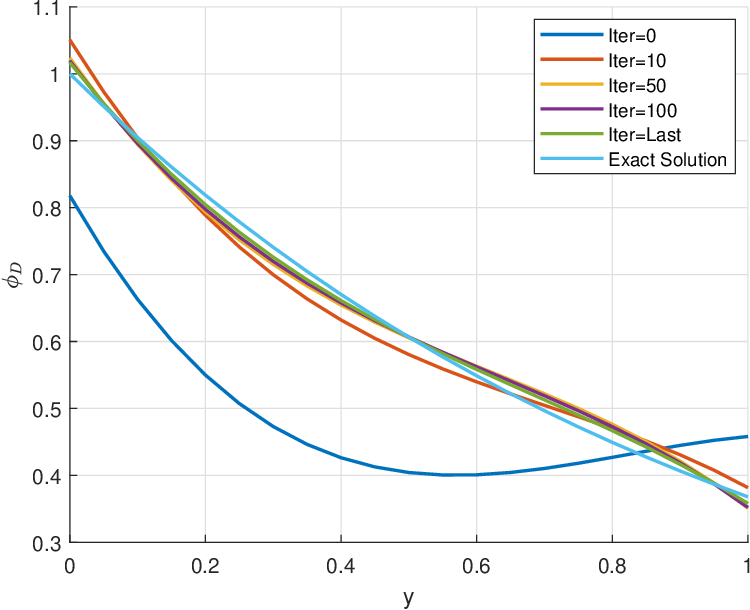}
        \caption{The trace of $u$ on $\Gamma_i$.}
        \label{fig2}
    \end{subfigure}
    \caption{The objective function $\mathcal{J_{DR}}(\phi_D)$ and the trace of $u$ on $\Gamma_i$ for various iterations with the best regularization parameter $\alpha=1e-08$.}
    \label{fig3}
\end{figure}

\begin{figure}[H]
    \centering
    \begin{subfigure}[b]{0.45\textwidth}
        \includegraphics[width=\textwidth]{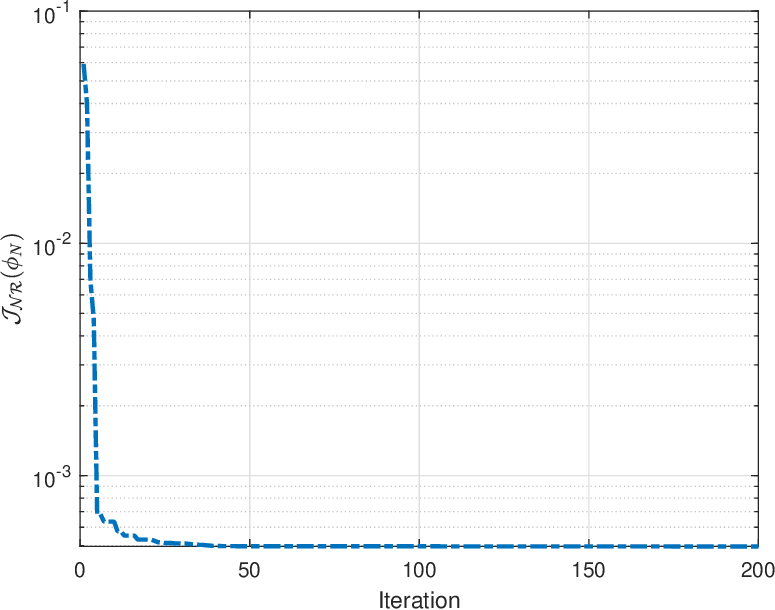}
        \caption{$\mathcal{J_{NR}}(\phi_N)$.}
        \label{fig4}
    \end{subfigure}
    \hfill
    \begin{subfigure}[b]{0.45\textwidth}
        \includegraphics[width=\textwidth]{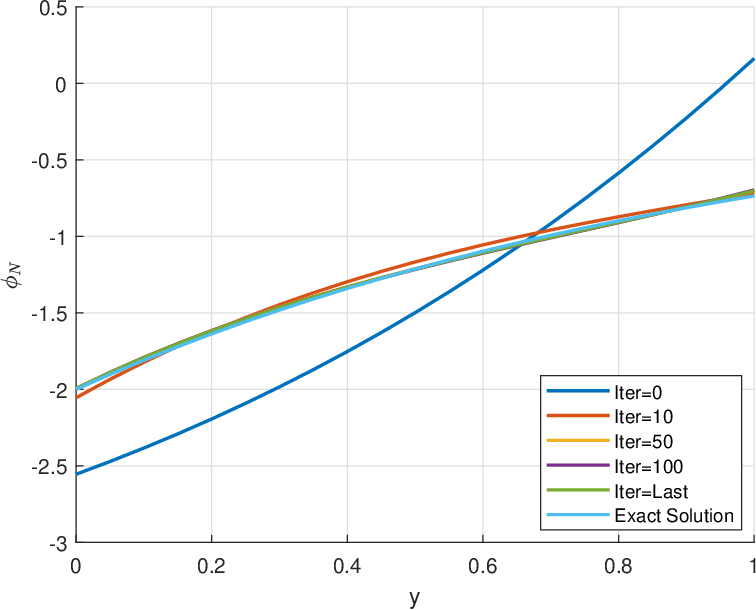}
        \caption{The derivative of $u$ $(\partial_n u)$ on $\Gamma_i$.}
        \label{fig5}
    \end{subfigure}
    \caption{The objective function $\mathcal{J_{NR}}(\phi_N)$ and the derivative of $u$ $(\partial_n u)$ on $\Gamma_i$ for various iterations with the best regularization parameter $\beta=1e-8$.}
    \label{fig6}
\end{figure}

Figures (\ref{fig2}) and (\ref{fig5}) illustrate the convergence of the numerical solution towards the analytical solution as the iterations progress. Initially, the numerical solution differs significantly from the exact solution, but the discrepancy diminishes rapidly with the number of iterations. This highlights the effectiveness of the iterative method in solving the inverse problem.

Additionally, Figures (\ref{fig1}) and (\ref{fig4}) show that the objective functions $\mathcal{J_{DR}}(\phi_D)$ and $\mathcal{J_{NR}}(\phi_N)$ decrease rapidly during the initial iterations. This suggests that the algorithm quickly approaches the minimum of the objective functions. As the iterations progress, the convergence rate slows down, but the objective functions still reach low values at $k=200$, ensuring a highly accurate solution with a good fit to the data.

\subsubsection{Stability of the Proposed Method}
\begin{figure}[H]
    \centering
    \begin{subfigure}[b]{0.45\textwidth}
        \includegraphics[width=\textwidth]{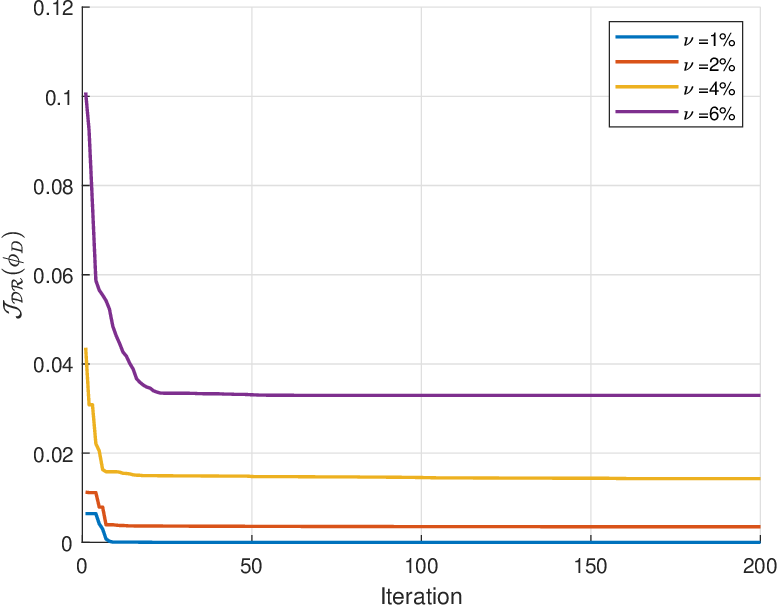}
        \caption{$\mathcal{J_{DR}}(\phi_D)$.}
        \label{fig7}
    \end{subfigure}
    \hfill
    \begin{subfigure}[b]{0.45\textwidth}
        \includegraphics[width=\textwidth]{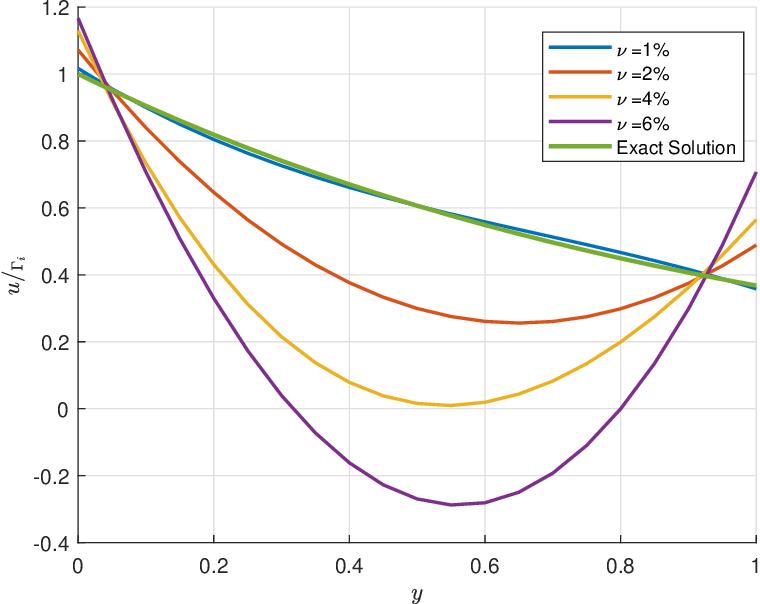}
        \caption{The trace of $u$ on $\Gamma_i$.}
        \label{fig8}
    \end{subfigure}
    \caption{The objective function $\mathcal{J_{DR}}(\phi_D)$ and the trace of $u$ on $\Gamma_i$ for various noise levels.}
    \label{fig9}
\end{figure}

\begin{figure}[H]
    \centering
    \begin{subfigure}[b]{0.45\textwidth}
        \includegraphics[width=\textwidth]{HelmHoltzCarre/DirichletObjectiveFunctionPSOTBruitCarreHelm.eps}
        \caption{$\mathcal{J_{NR}}(\phi_N)$.}
        \label{fig10}
    \end{subfigure}
    \hfill
    \begin{subfigure}[b]{0.45\textwidth}
        \includegraphics[width=\textwidth]{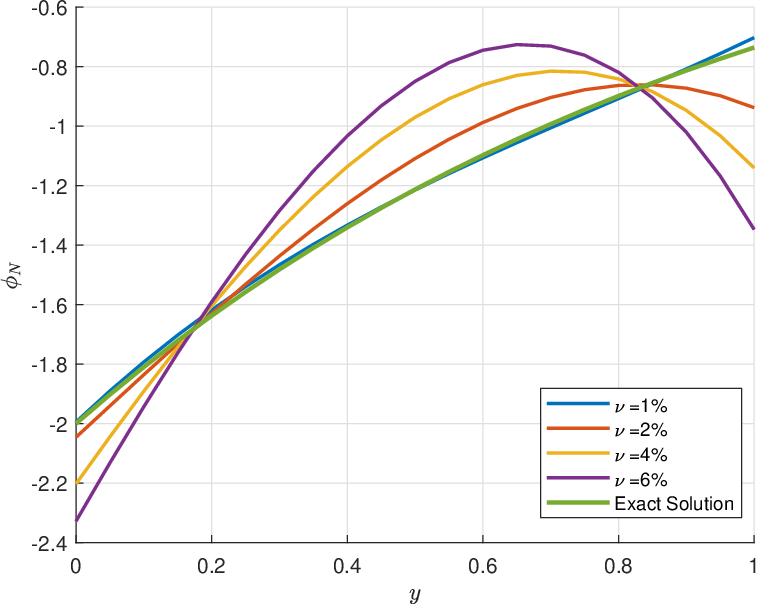}
        \caption{The derivative of $u$ $(\partial_n u)$ on $\Gamma_i$.}
        \label{fig11}
    \end{subfigure}
    \caption{The objective function $\mathcal{J_{NR}}(\phi_N)$ and the derivative of $u$ $(\partial_n u)$ on $\Gamma_i$ for various noise levels.}
    \label{fig12}
\end{figure}

Figures (\ref{fig9}) and (\ref{fig12}) provide insights into the robustness of the numerical approximations against varying noise levels in the measurement data. The results indicate that the numerical solutions obtained through the proposed method remain consistent with the analytical solution, even under different noise conditions.

Figures (\ref{fig8}) and (\ref{fig11}) compare the numerical and analytical solutions under different noise levels. While the numerical solution deviates slightly from the exact solution as noise increases, the discrepancies remain minimal even at noise levels of $3\%$. 

Figures (\ref{fig7}) and (\ref{fig10}) illustrate the cost functions at noise levels $\nu = 0\%, 1\%, 2\%, 3\%$. As noise increases, the cost function rises, indicating a less precise fit to the data. However, the cost functions remain relatively low, suggesting that the numerical solution maintains a good fit even under significant noise.

\subsection{Unit Disc}
\subsubsection{Regularization Parameter Choice}

\begin{table}[h!] 
\centering
\begin{tabular}{|c|c|c|c|c|c|c|}
\hline
$\eta_{\alpha,\beta}$        & 1e-03       & 1e-04       & 1e-05       & {\color{red} 1e-06}       & 1e-07       & {\color{green} 1e-08}       \\ \hline
$\mathcal{J_{DR}}(\phi_D)$     & 9.2808e-05 & 1.2965e-03  & 3.7476e-04  & 8.9112e-06  & 7.5673e-04  & {\color{green} 6.7854e-06} \\ \hline
$\mathcal{J_{NR}}(\phi_N)$ & 2.1708e-03  & 2.1453e-03  & 2.1456e-03  & {\color{red} 2.143e-03}  & 2.1577e-03  & 2.1468e-03 \\ \hline
\end{tabular}
\caption{Values of the cost functions $\mathcal{J_{DR}}(\phi_D)$ and $\mathcal{J_{NR}}(\phi_N)$ for different regularization parameters $\eta$.}
\label{RegCercle}
\end{table}

Table \ref{RegCercle} presents the values of the cost functions $\mathcal{J}_{DR}(\phi_D)$ and $\mathcal{J}_{NR}(\phi_N)$ for different regularization parameters $\eta_{\alpha,\beta}$, illustrating the impact of this parameter on the quality of the reconstruction. In the Dirichlet case, the cost function exhibits clear sensitivity to $\eta$, with the lowest value, $\mathcal{J}_{DR}(\phi_D) = 6.7854 \times 10^{-6}$, achieved at $\eta = 10^{-8}$, indicating the most accurate reconstruction. In contrast, a significantly higher value of $1.2965 \times 10^{-3}$ is observed at $\eta = 10^{-4}$, highlighting a non-monotonic trend and the potential degradation of reconstruction quality when $\eta$ is not appropriately chosen. For the Neumann formulation, the cost function $\mathcal{J}_{NR}(\phi_N)$ remains relatively stable across all tested values, with the best result, $2.143 \times 10^{-3}$, obtained at $\eta = 10^{-6}$. These results show that the Dirichlet approach is more sensitive to the choice of the regularization parameter, whereas the Neumann formulation demonstrates greater robustness. Based on these observations, the optimal values $\eta_{\alpha,\beta} = 10^{-8}$ for the Dirichlet case and $\eta_{\alpha,\beta} = 10^{-6}$ for the Neumann case are selected and used throughout the remainder of this study to ensure accurate and stable reconstructions.

\subsubsection{Convergence of the Algorithm}
\begin{figure}[H]
    \centering
    \begin{subfigure}[b]{0.4\textwidth}
        \includegraphics[width=\textwidth]{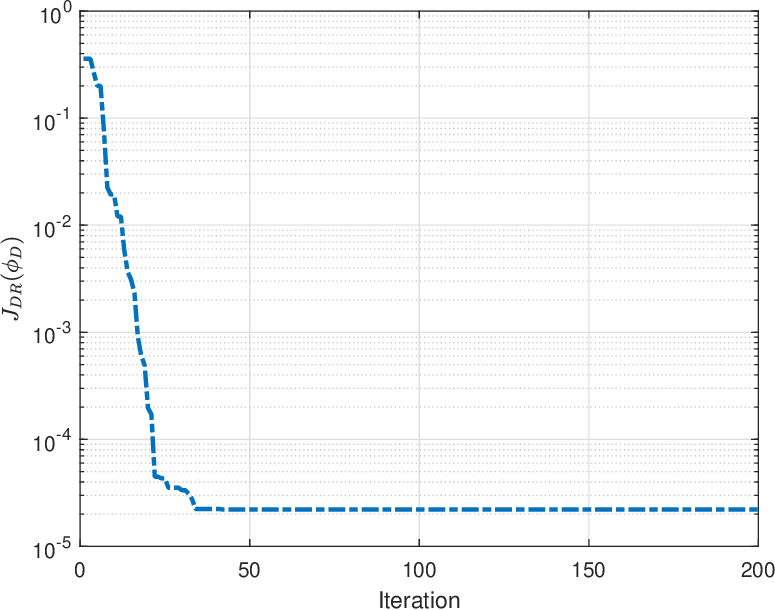}
        \caption{$\mathcal{J_{DR}}(\phi_D)$}
        \label{fig13}
    \end{subfigure}
    \hfill
    \begin{subfigure}[b]{0.4\textwidth}
        \includegraphics[width=\textwidth]{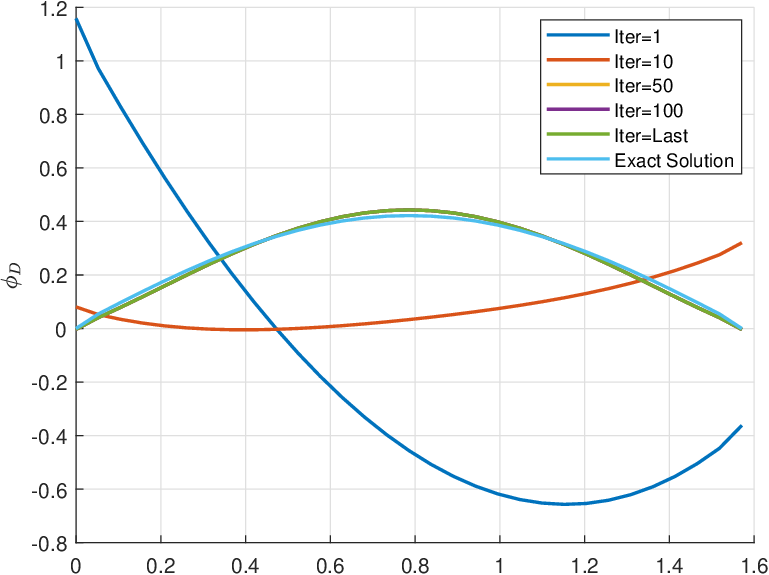}
        \caption{The trace of $u$ on $\Gamma_i$.}
        \label{fig14}
    \end{subfigure}
    \caption{The objective function $\mathcal{J_{DR}}(\phi_D)$ and the trace of $u$ on $\Gamma_i$ for various iterations with the best regularization parameter $\alpha=1e-08$.}
    \label{fig15}
\end{figure}

\begin{figure}[H]
    \centering
    \begin{subfigure}[b]{0.4\textwidth}
        \includegraphics[width=\textwidth]{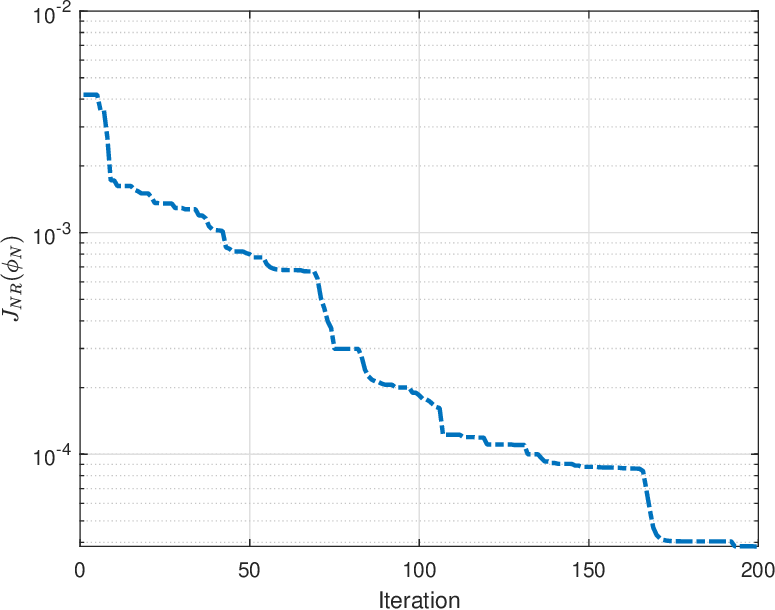}
        \caption{$\mathcal{J_{NR}}(\phi_N)$}
        \label{fig16}
    \end{subfigure}
    \hfill
    \begin{subfigure}[b]{0.4\textwidth}
        \includegraphics[width=\textwidth]{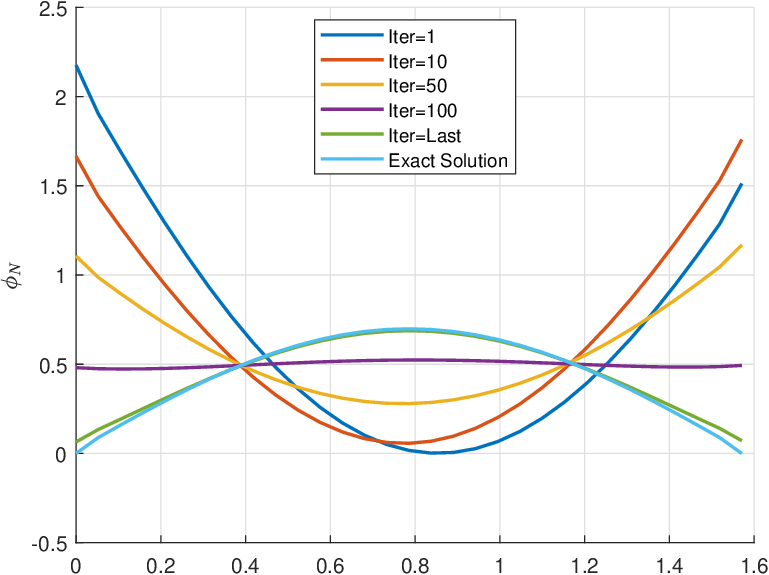}
        \caption{The derivative of $u$ $(\partial_n u)$ on $\Gamma_i$.}
        \label{fig17}
    \end{subfigure}
    \caption{The objective function $\mathcal{J_{NR}}(\phi_N)$ and the derivative of $u$ $(\partial_n u)$ on $\Gamma_i$ for various iterations with the best regularization parameter $\beta=1e-06$.}
    \label{fig18}
\end{figure}

Figures (\ref{fig14}) and (\ref{fig17}) show that the numerical solution converges to the exact solution as the iteration progresses. Initially, the numerical solution differs significantly from the exact solution, but the difference decreases rapidly as the iterations increase. This demonstrates the effectiveness of the iterative method in solving the inverse problem.

Additionally, Figures (\ref{fig13}) and (\ref{fig16}) illustrate that the objective functions $\mathcal{J_{DR}}(\phi_D)$ and $\mathcal{J_{NR}}(\phi_N)$ decrease rapidly during the initial iterations, indicating that the algorithm is progressing toward the minimum of the objective functions. As the iterations continue, the convergence slows, but the objective function eventually reaches a low value at $k=200$, signifying a highly accurate solution with a good fit to the data.

\subsubsection{Stability of the Proposed Method}
\begin{figure}[H]
    \centering
    \begin{subfigure}[b]{0.4\textwidth}
        \includegraphics[width=\textwidth]{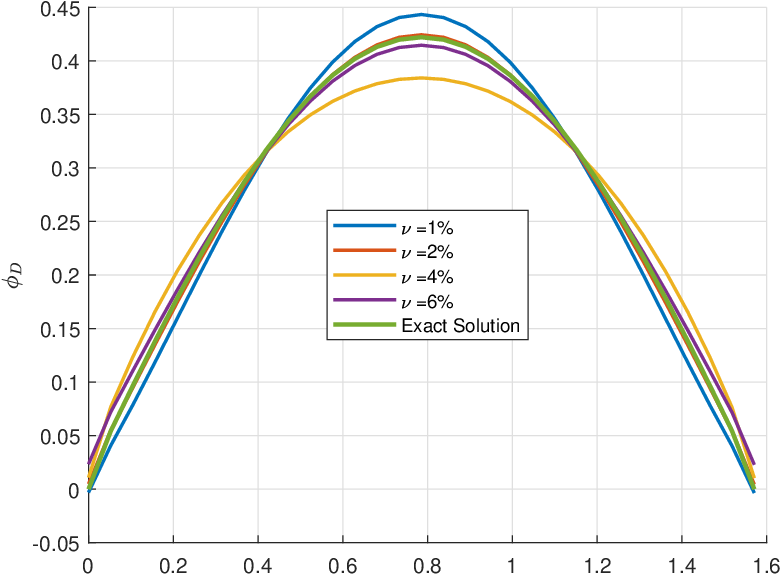}
        \caption{$\mathcal{J_{DR}}(\phi_D)$}
        \label{fig19}
    \end{subfigure}
    \hfill
    \begin{subfigure}[b]{0.4\textwidth}
        \includegraphics[width=\textwidth]{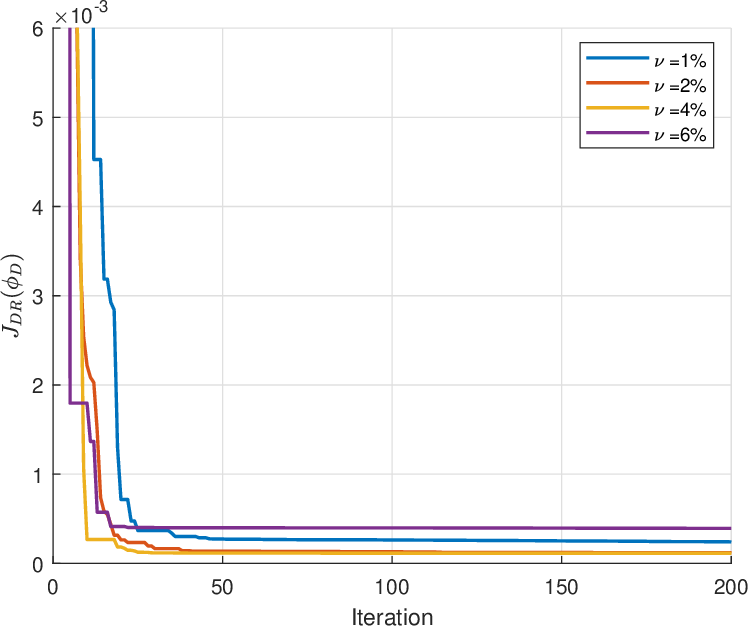}
        \caption{The trace of $u$ on $\Gamma_i$.}
        \label{fig20}
    \end{subfigure}
    \caption{The objective function $\mathcal{J_{DR}}(\phi_D)$ and the trace of $u$ on $\Gamma_i$ for various levels of noise.}
    \label{fig21}
\end{figure}

\begin{figure}[H]
    \centering
    \begin{subfigure}[b]{0.4\textwidth}
        \includegraphics[width=\textwidth]{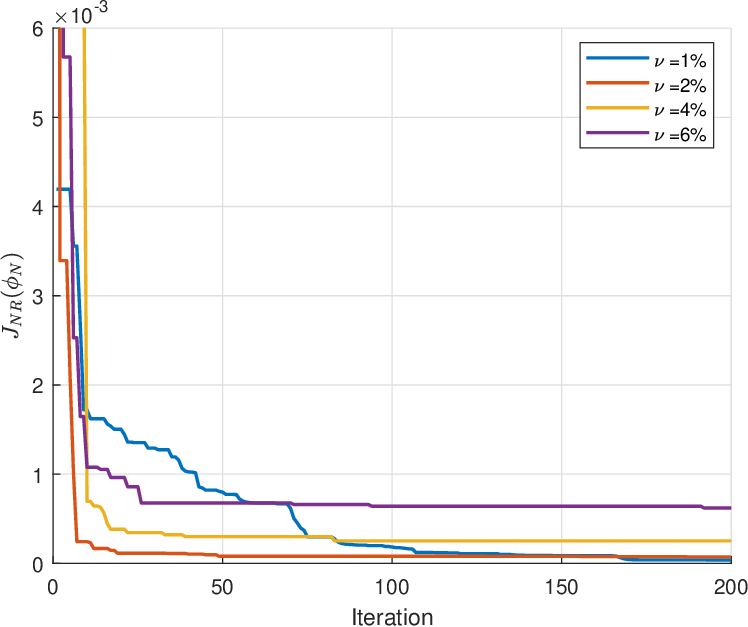}
        \caption{$\mathcal{J_{NR}}(\phi_N)$}
        \label{fig22}
    \end{subfigure}
    \hfill
    \begin{subfigure}[b]{0.4\textwidth}
        \includegraphics[width=\textwidth]{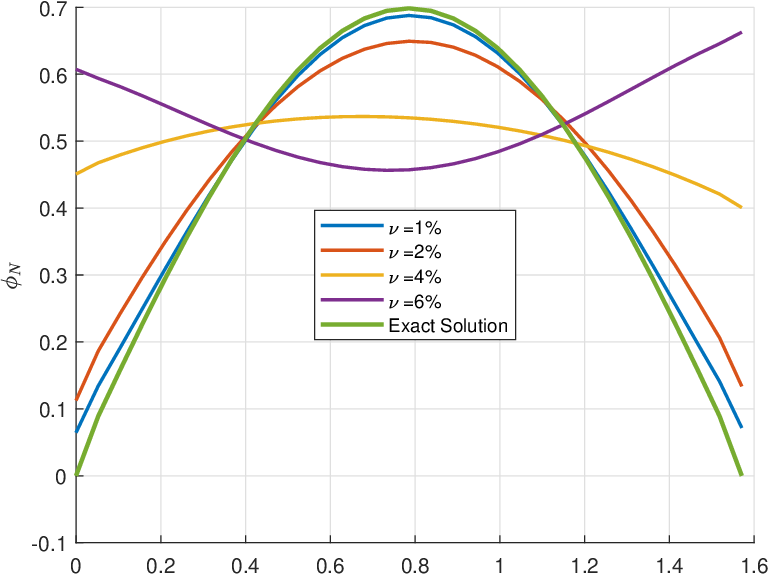}
        \caption{The derivative of $u$ $(\partial_n u)$ on $\Gamma_i$.}
        \label{fig23}
    \end{subfigure}
    \caption{The objective function $\mathcal{J_{NR}}(\phi_N)$ and the derivative of $u$ $(\partial_n u)$ on $\Gamma_i$ for various levels of noise.}
    \label{fig24}
\end{figure}

The results presented in Figures (\ref{fig21}) and (\ref{fig24}) shed light on the robustness of the numerical approximations to varying noise levels in the measurement data. As shown, the numerical solutions remain consistent with the exact solution, even when different levels of noise are introduced into the measurement data.

Figures (\ref{fig20}) and (\ref{fig23}) compare the numerical solution to the exact solution for varying noise levels in the measurement data. While the numerical solution deviates slightly from the exact solution as the noise increases, the discrepancy remains small even when the noise level reaches $3\%$.

Figures (\ref{fig19}) and (\ref{fig22}) display the cost function under various noise levels. As the noise level increases, the cost function increases, suggesting a less accurate fit to the data. However, even at higher noise levels, the cost function remains relatively low, indicating that the numerical solution continues to provide a good fit to the data.

It is important to note that reconstructing the Dirichlet condition on the inaccessible part of the boundary is typically more accurate and reliable than reconstructing the Neumann condition in both cases. This difference in accuracy can depend on various factors, including the geometry of the domain or the regularity of the solution being reconstructed.

\section{Conclusion}
\label{sec:5}
In this paper, we focus on a challenging ill-posed inverse problem, namely the data completion problem for the Helmholtz equation. While several optimization methods are available for approximating the solution of such problems, this study investigates the potential of Particle Swarm Optimization (PSO), which does not require assumptions about the regularity of the considered functional. We propose an optimization formulation that incorporates a Tikhonov regularization term and demonstrate the effectiveness of the proposed approach through numerical simulations on both regular and irregular domains. The results indicate that the PSO method can successfully address the problem, even in situations involving domain irregularities.

Our study contributes to the growing body of research on solving ill-posed inverse problems using optimization techniques. The PSO-based approach offers an alternative to existing methods that require assumptions about the regularity of the functional. Future research could explore the combination of PSO with other optimization techniques to improve the accuracy and efficiency of solving ill-posed inverse problems.

\end{document}